\documentclass[12pt]{article}

\usepackage{amscd}
\usepackage{amssymb,amsmath}

\newtheorem{theorem}{Theorem}
\newtheorem{prop}[theorem]{Proposition}
\newtheorem{lemma}[theorem]{Lemma}
\newtheorem{cor}[theorem]{Corollary}

\date{Sept 30, 2008}
\author{Vladimir Baranovsky}
\title{Bundles on non-proper schemes: representability}

\begin{document}

\maketitle

\section{Introduction}

Let $X_k$ be a proper irreducible separated scheme of finite type over a field
$k$. We will also assume that $X_k$ satisfies Serre's $S_2$ condition 
(for the sake of simplicity the reader may think that $X_k$ is smooth or a
 locally complete intersection).
For a noetherian scheme $T$ over $k$ denote $X_T = (X_k) \times_{Spec(k)} T$  and let $\Phi(T)$ be the collection of all closed
subsets $Z\subset X_T$ such that every point $z \in Z$ has codimension
$\geq 3$ in its fiber over $T$.  Fix a reductive group $G$ over $k$.

\bigskip
\noindent 
\textbf{Definition.} In the notation above, let
$F_G(T)$ be a groupoid category with the objects $(E, U)$ where 
$E$ is a principal $G$-bundle defined on an open subscheme 
$U\subset X_T$, such that the closed complement of $U$ is a 
subset in $\Phi(T)$. 
A morphism $(E_1, U_1) \to (E_2, U_2)$ is an isomorphism $E_1|_W\simeq E_2|_W$
on an open subset $W \subset U_1 \cap U_2$ such that the complement of $W$
is again in $\Phi(T)$. The composition of morphisms is defined in an obvious way.

\bigskip
\noindent
For any morphism $\alpha: T' \to T$ of noetherian schemes over $k$ we have 
pullback functors $\alpha^*: F_G (T)\to F_G(T')$ satisfying the usual 
compatibility conditions for 
any pair of morphisms $T'' \stackrel{\beta}\longrightarrow T'
\stackrel{\alpha} \longrightarrow T$, i.e. $F_G$ is a groupoid
over the category of noetherian $k$-schemes, cf. Section 1 in 
\cite{Ar}. As usual, we 
will mostly deal with its restriction to affine noetherian 
$k$-schemes,
 writing $X_A$, $F(A)$ and $\Phi(A)$ instead of $X_{Spec(A)}$,
$F(Spec(A))$ 
and $\Phi(Spec(A))$,  respectively.
The main goal of this paper is the
following result.

\begin{theorem}
$F_G$ is an algebraic stack, locally of finite type over $k$, 
with separated  and quasi-compact diagonal.
\end{theorem}
Thus we obtain a partial ``compactification" of the stack of 
$G$-bundles on $X_k$ (although we are not claiming here that the 
stack of $G$-bundles in dense in $F_G$).
Our strategy of proof is straightforward, if seldom used: we apply Artin's
representability criterion, cf. Theorem 5.3 in \cite{Ar} for a statement and
\cite{L}, \cite{Ao} for examples of application.

For the most part of the paper (see Sections 3-7) we consider the case
of vector bundles, i.e. work with $G = GL(r)$ for fixed $r \geq 1$; and 
write $F$ instead of $F_G$. In Section 8 we show 
how the proof is extended to the case of general $G$ and also explain why
the result fails when the ``codimension 3" condition in the definition
of $\Phi(T)$ is replaced by ``codimension 2".

\medskip
\noindent
\textbf{Acknowledgements.} The main result of this paper was conjectured 
by Vladimir Drinfeld whom the author thanks for the useful conversations and
encouragement. Thanks are also due to Victor Ginzburg for helping to clarify 
the role of Serre's $S_2$ condition (which is closely related to 
the concept of a reflexive sheaf) 
and to Max Lieblich for his comments on stacks.  
This work was supported by the Sloan
Research Fellowship and UCI Teaching Relief Grant.

\section{Depth and local cohomology.}

We use definitions and basic properties of local cohomology
and depth which can be found e.g. in Section 18 and Appendix 4 to \cite{E} and
Sections I-III of [SGA\;2].  Observe that the scheme $X_T$
will not satisfy Serre's $S_2$ condition $depth_x \mathcal{O}_{X_T}
\geq min(2, \dim \mathcal{O}_{X_T, x})$, since no depth assumptions 
are imposed on $T$. However, we can formulate a relative version 
of this condition.

\bigskip
\noindent
\textbf{Definition.} For any point $x \in X_T$ let $d(x)$ be the codimension
of $x$ in its fiber over $T$. 

\begin{lemma}
\label{serre}
Let $E$ be a vector bundle on an open subset $U \subset X_T$ and
$M$ a coherent sheaf on $T$. For $t \in T$ let 
$X_t$ be the fiber of $X_T$ over $t$ and suppose that for every $t$ the 
structure sheaf of the intersection
$U_t = U\cap X_t$ satisfies Serre's condition $S_n$. Set 
$E_M= E \otimes_{\mathcal{O}_T} M$. Then for any $x \in U$ one has 
$depth_{U, x} E_M \geq min(n, d(x))$.

\end{lemma}
\textbf{Definition.} 
We will call the inequality stated in this lemma \textit{the relative
$S_n$ condition}. In this paper $n = 2$ or $3$.

\medskip
\noindent
\textit{Proof of the lemma.} Since the question is local we can 
 take $E = \mathcal{O}_U$. 
Thus we assume $T= Spec(A)$, 
$U = X = Spec(C) \times_{Spec(k)} Spec(A)$ and 
$E_M$ of the form $C \otimes_k M$. Since the fiber over $t \in T$ is 
given by $X_{k(t)}$, the codimension of $x$ in $X_t$ is equal to 
the codimension of its image $x' \in X_k$. Setting $r = min(n, d(x))$
we can find a $C$-regular sequence $f_1, \ldots, f_r$ in the maximal
ideal of $x'$ in $C$. The same $f_i$ viewed as elements of 
$C \otimes_k A$ will belong to the maximal ideal of $x$ and form 
a $C \otimes_k M$-regular sequence, which finishes the proof. $\square$

\bigskip
\noindent
Below, dealing with obstructions, deformations and infinitesimal
automorphisms we need the following construction. For for any coherent sheaf 
$E$ on $X_T$ set
$$
H^i_{T, \Phi} (E) = \underrightarrow{\lim}_{Z \in \Phi(T)} 
H^i(X_T \setminus Z, E)
$$
where the filtered direct limit is taken with respect to the inclusion
of closed subsets $Z \subset Z'$ in $\Phi(T)$.
If $T = Spec(A)$ is affine, we write $H^i_{A, \Phi}(E)$ instead 
of $H^i_{Spec(A), \Phi}(E)$. If $E$ is defined only on an open subset
$U = X_T \setminus Z_0$ with $Z_0\in \Phi(T)$ we can use the same definition
but take the limit over those $Z$ which contain $Z_0$.
Observe that for $i =0, 1$ the cohomology groups in the limit in fact
stabilize under certain restrictions on $E$:

\begin{lemma}
\label{stabilize}
With the notation just introduced, assume that $Z \subset Z'$ are
in $\Phi(T)$. If $E$ satisfies the relative $S_2$ condition 
on $X_T \setminus Z$ then the natural restriction morphism
$$
\rho_i: H^i(X_T \setminus Z, E) \to H^i(X_T \setminus Z', E) 
$$
is an isomorphism for $i =0$ and injective for $i = 1$. If in addition
$E$ satisfies the relative $S_3$ condition on $X_T \setminus Z$
 then $\rho_i$ is an isomorphism for $i = 0, 1$
and injective for $i = 2$. 
\end{lemma}
\textit{Proof.} Denote $U = X_T \setminus Z$, $W = U \cap Z'$ and 
consider the spectral sequence of local cohomology
$H^p(U, \mathcal{H}^q_W(E)) \Rightarrow H^{p+q}_W(U,E)$. By the 
relative $S_2$ condition the local cohomology sheaves $\mathcal{H}^i_W(E)$
vanish for $i=0,1$ while the
relative $S_3$ condition at the points of $W$ also implies
$\mathcal{H}^2_W(E) = 0$. 
Now the assertion follows from the standard long exact sequence
$$
\ldots \to H^i_{W}(U, E) \to H^i(U, E) \to H^i(U \setminus W, E) 
\to H^{i+1}_W (U, E) \to \ldots \hspace{2cm} \square
$$

\bigskip
\noindent
Observe that the relative $S_3$ condition on $E_M$ also holds if $Z$ 
contains 
$$
Z^\circ_T = Z^\circ \times_{Spec(k)} T \in \Phi(T)
$$
where $Z^\circ \subset X_k$ is the set of all points 
 in $X_k$ where the $S_3$ condition fails for the structure sheaf. Observe that
 $Z^\circ$ is closed by [EGA\;$\textrm{IV}_2$], Proposition 6.11.2; 
since its complement contains all points
of codimension $\geq 2$ we indeed have $Z^\circ_T \in \Phi(T)$.

\begin{cor}
\label{finitely}
If $Z^\circ_T \subset Z$ then 
$$
H^i_{T, \Phi} (E_M) = H^i(X_T \setminus Z, E_M)
$$
for $i = 0, 1$ and $E_M= E\otimes_{\mathcal{O}_T} M$ as in Lemma
\ref{serre}. Moreover, if $j$ 
stands for the open embedding $X_T \setminus Z \hookrightarrow X_T$ then the 
sheaves $j_*(E_M)$ and $R^1j_*(E_M)$ are coherent on $X_T$ and 
if $T = Spec(A)$ is affine the two stable cohomology groups are finitely 
generated $A$-modules.
\end{cor}
\textit{Proof.} Stabilization follows immediately from the previous lemma.
Coherence of the two direct images is due to [SGA 2], VII.2.3 
while the finite generation is proved by combining the
spectral sequence 
$H^p(X_T, R^qj_*(E_M)) \Rightarrow H^{p+q}(X_T \setminus Z, E_M)$
with the fact that $X_T$ is proper over $T$. $\square$

\bigskip
\noindent
\textbf{Remark.}
Of course, for $i=2$
even the cohomology group $H^2(\mathbb{P}_k^3 \setminus P, \mathcal{O})$
is infinite dimensional over $k$ for any closed point $P$.

\section{Locally finite presentation.}

In this section we do not use the $S_2$ assumption on $X_k$.
Let $R = \underrightarrow{\lim} \;R_\alpha$ be a 
filtered direct limit of $k$-algebras. 
\begin{prop}
$$
\underrightarrow{\lim}\; F(R_\alpha) \to F(R)
$$
is an equivalence of categories. 
\end{prop}
\textit{Proof.} The assertion
means that any object $(E, U) \in F(R)$ is an image
of some $(E_\alpha, U_\alpha)\in F(R_\alpha)$ and that, 
whenever $(E_\alpha, U_\alpha)$ and $(E_\beta, U_\beta)$ give isomorphic
objects in $F(R)$, there exists $\gamma$ such that $\gamma \geq \alpha$, 
$\gamma \geq \beta$ and the corresponding objects in $F(R_\gamma)$ are 
isomorphic. In addition, a similar condition should hold for morphisms.

To prove the assertion for objects, consider a vector bundle $E$ on
$U \subset X_R$ and take a finite affine covering
$\{U_i\}$ of $U$ such that $E|_{U_i}$ is trivial. 

Using the results of Sections 8.2-8.5 of [EGA $\textrm{IV}_3$] we see that 
there exists
$\alpha$ and open subsets $U^\alpha_i$ such that $U_i = \pi_{\alpha}^{-1}
(U_i^{\alpha})$ where $\pi_\alpha: X_R \to X_{R_\alpha}$
is the natural projection. Since in general a scheme $W$ is affine iff the 
canonical morphism $W\to Spec(\Gamma(W, \mathcal{O}_W))$ is an isomorphism,
by increasing $\alpha$ if necessary we can assume that all $U^\alpha_i$ are
affine. 

The transition functions for $E$ given by $\phi_{ij}:U_i \cap U_j \to GL_r(k)$
can be viewed as automorphisms of the trivial bundle. Increasing
$\alpha$ we can assume that $\phi_{ij}$ arise from 
regular maps $\phi^\alpha_{ij}: U^\alpha_i \cap U^\alpha_j \to GL(r)$. 
Increasing $\alpha$ again we can assume that $\phi^\alpha_{ij}$ 
satisfy the cocycle condition and thus define a vector bundle $E^\alpha$
on $U^\alpha = \bigcup_{i} U^\alpha_i$. By construction, $E \simeq \phi_\alpha^*
E^\alpha$. To show that the closed complement $Z^\alpha$ of $U^\alpha$ 
is in $\Phi(Spec(R_\alpha))$ (again, after 
a possible increase of $\alpha$)  note that 
$U = \pi_{\alpha}^{-1} (U^\alpha)$ and the fibers of $X_R
\to Spec(R)$ are obtained from the fibers 
$X_{R_\alpha} \to Spec(R_\alpha)$ by extension of scalars. Therefore
the closed subset $W$ of points $s \in Spec(R_\alpha)$ for which codimension 
of $Z^\alpha \cap X_s$  is $\leq 2$ has empty preimage in 
$Spec(R)$. Therefore, for some $\alpha' \geq \alpha$ the preimage of
$W$ in $Spec(R_{\alpha'})$ is empty and we can replace $\alpha$ by $\alpha'$. 

To prove surjectivity on morphisms, let $(E_1, U_1)$ and $(E_2, U_2)$ 
be two objects in $F(R)$ and 
suppose we are given an isomorphism $\phi: E_1|_U \simeq
E_2|_U$ where $U \subset U_1 \cap U_2$ is open 
with its closed complement in $\Phi(R)$. 

By the previous argument, we can assume that $E_i$ is isomorphic to the 
pullback of some vector bundle $E^\alpha_i$ on an open subset $U^\alpha_i
\subset X\times Spec(R_\alpha)$. Increasing $\alpha$ we can assume that 
$U$ is the preimage of an open subset $U^\alpha \subset U^\alpha_1 \cap
U^\alpha_2$. Then $E^\alpha_1|_{U^\alpha}$ and $E^\alpha_2|_{U^\alpha}$ 
become isomorphic after pullback to $U$ hence by \textit{loc. cit.} 
by increasing $\alpha$ we can find an isomorphism $\phi^\alpha:
E^\alpha_1|_{U^\alpha} \simeq E^\alpha_2|_{U^\alpha}$ which induces
$\phi$ on $U$. As before, we may have to increase $\alpha$ one more time
to ensure that the complement of $U^\alpha$ 
is in $\Phi(R_\alpha)$. 

Injectivity on morphisms is an immediate consequence of Theorem 8.5.2 in 
\textit{loc. cit}. $\square$

\section{Small affine pushouts.}

Let $A_0$ be a noetherian $k$-algebra, and $A' \to A$ a surjection of two
infinitesimal extensions of $A_0$ such that $M= ker(A' \to A)$ is a
finite $A_0$ module. Let $B$ be a noetherian ring and $B\to A$
a morphism, such that the composition $B \to A \to A_0$ is surjective. 

Denote by $B'$ the pushout $A' \times_A B$, i.e. the subset of pairs
$(a, b) \in A' \times B$ which have the same image in $A$. Then $B' \to B$
is surjective and its kernel may be identified with $M$ viewed as 
a $B$-module. Observe that $Spec(B')$ is homeomorphic to $Spec(B)$, while
$Spec(A'), Spec(A)$ and $Spec(A_0)$ are homeomorphic to each other, and 
$Spec(A_0) \to Spec(B)$ is naturally a closed subscheme by the assumption.

\medskip
\noindent
Fix an object $a = (E, U(A)) \in F(A)$. Let $F_a(B)$ the 
groupoid of extensions of $a$ over $Spec(B)$, and similarly for $A', 
B'$.
\begin{prop}
The natural functor
$$
F_a (B') \to F_a(A') \times F_a(B)
$$
is an equivalence of groupoids. 
\end{prop}
\textit{Proof.}
Suppose that $(E', U(A'))$, $(E'', U(B))$ are two extensions of 
$(E, U(A))$ to $X_{A'}$ and $X_B$, respectively. Since $Spec(A)$ and $Spec(A')$ are homeomorphic, 
shrinking $U(A'), U(B)$ and $U(A)$, if necessary, we can assume 
$U(A') \simeq U(A) \simeq U(B) \cap 
X_{A_0}$ (homeomorphisms induced by the natural embeddings). 
Denote by $U(B')$ the subset $U_B$ viewed as an open subset of
$X_{B'}$. We have a commutative diagram
$$
\begin{CD}
U(A) & @>p>> & U(A') \\  
@V{q}VV && @V{i}VV \\
U(B) & @>j>> & U(B')
\end{CD}
$$
where the horizontal arrows are homeomorphisms. 
Since $q^* E'' \simeq E$, $p^* E'\simeq E$ and 
$i p = j q$ there will be an exact sequence on $U(B')$
$$
i_* E' \oplus j_* E'' \to (ip)_* E \to 0
$$ 
where the first arrow is given by the difference of the obvious 
canonical maps. One can check that the
kernel $E'''$ of the first arrow is a locally free sheaf of rank $r$ on 
$U(B') \subset X_{B'}$ such that $i^* E'''\simeq E'$, 
$j^* E''' \simeq E''$, 
$(ip)^* E'''\simeq E$ in a compatible way. A further 
straightforward check shows
that the correspondence $(E', E'') \mapsto E'''$ induces 
which is an equivalence of categories. $\square$

\section{Automorphisms, deformations, obstructions.}

As in Section 2, for a vector bundle $E$ on an 
open subset of $X_A$ or $X_{A_0}$ and an $A_0$-module
$M$ we will write $E_M$ for the tensor product of $E$
with the pullback of $M$ viewed as a coherent sheaf on 
$Spec(A)$ or $Spec(A_0)$, respectively. 

\subsection{$Aut$, $D$, $O$}

\textit{Infinitesimal automorphisms.}

\medskip
\noindent
Let $A = A_0$, $A' = A_0 \oplus M$.
Any $a_0 \in F(A_0)$ given by a pair $(E, U)$ admits a trivial extension to 
$A_0 \oplus M$ defined by $E'= E \oplus E_M$. We are 
interested in the group of automorphisms $Aut_{a_0} (A_0 + M)$ of the
bundle $E'$, which restrict to identity over $A_0$. Every such automorphism
is defined uniquely by a morphism $E \to E_M$ given, perhaps,
on a smaller open subset $V \subset U$. In other words
$$
Aut_{a_0} (A_0 + M)= H^0_{A_0, \Phi}(End(E)_M)
$$
By Corollary \ref{finitely} this is a finitely generated module over $A_0$.

\bigskip
\noindent
\textit{Deformations}

\medskip
\noindent
Now consider $D_{a_0}(M)$,  the set of isomorphism classes
of extensions of $a_0 = (E, U)$ to $A'$.  A standard argument, cf. e.g.
Chapter IV of \cite{I}, identifies $D_{a_0}(M)$ with 
$H^1_{A_0, \Phi} (End(E)_M)$ which is 
also finitely generated over $A_0$, as established in Corollary \ref{finitely}.

\bigskip
\noindent
\textit{Obstructions}

\medskip
\noindent
Our goal is to define a finitely generated submodule
of $H^2_{A_0, \Phi}(End(E)_M)$ which will serve as
obstruction module for our problem.
For a general square zero extension  
$0 \to M \to A' \to A \to 0$ and
 $a = (E, U) \in F(A)$, let 
$U' \subset X_{A'}$ be the open subset homeomorphic to $U$, 
with its natural structure of an open subscheme of $X_{A'}$.
By [I] there is an obstruction to deforming $E$ over $U'$, given by 
a class 
$$
\omega (E) \in Ext^2_{U_A}(End(E), M)
$$
which is a Yoneda product of two classes
$$
a(E) \in Ext^1_{U}(End(E),  L_{U}),  \qquad \kappa(U/U')  
\in Ext^1_{U} (L_{U}, M).
$$
Here $L_{U}$ is the cotangent complex of $U$ over $k$,
 $a(E)$ is the Atiyah class of $E$ and 
$\kappa(U/U') \in Ext^1_{U} (L_{U}, M) $ is 
the Kodaira-Spencer class of $U'$ viewed as a deformation of $U$, 
cf. \textit{loc.cit.} 
Moreover, since $X_A$ is a direct product of $X_k$ and $Spec(A)$, its cotangent
 complex over $k$ 
splits into a direct sum $L_{X_k} \oplus L(A)$ of the pullbacks of 
cotangent complexes from $X_k$ and $Spec(A)$, respectively. Since $X_{A'}$
and $U'$ viewed as deformations of $X_A$ and $U$, respectively, are 
induced by $Spec(A) \to Spec(A')$, the Kodaira-Spencer class 
$\kappa(U/U')$ belongs to the direct factor 
$Ext^1_{U}(L(A), M) \subset Ext^1_U(L_{U},
M)$. Therefore, $\omega(E)$ may be viewed as the product of the 
Kodaira-Spencer class with the image of the Atiyah class
$a'(E) \in Ext^1_{U}(End(E), L(A))$. 
If we represent $a'(E)$ by an extension
$$
0 \to L(A) \to Q \to End(E) \to 0,
$$
then $\omega(E)$ will become the image of $\kappa(U/U')$ under the 
connecting homomorphism:
$$
\ldots \to Ext^1_{U}(Q, M) \to Ext^1_{U}(L(A), M) \to 
Ext^2_{U}(End(E), M) \to \ldots
$$
Therefore, the obstruction $\omega(E)$ can be viewed as a element
of the following cokernel $O_a(U, M)$:
$$
Ext^1_{U}(Q, M) \to Ext^1_{U}(L(A), M) \to O_a(U, M) \to 0.
$$
Since the cotangent complex of $Spec(A)$ is concentrated in non-positive
degrees and has finitely generated cohomology, we can find a 
free resolution $\ldots \to L_2 \to L_1 \to L_0 \to L(A) \to 0$ and therefore
a locally free resolution $\ldots \to L_2 \to L_1 \to Q_0 \to Q \to 0$.
The standard spectral sequence $E_1^{p, q} = Ext^p_{U}(L_q, M)
\Rightarrow Ext^{p+q}_{U}(L(A), M)$ shows that 
$Ext^1_{U}(L(A), M)$ has a
two-step filtration 
with associated graded depending only on 
$E_1^{p, q}$ for $p \leq 1$ and
$q \leq 2$ which by Corollary \ref{finitely}
implies that it is finitely generated over $A$.
Similar argument applies to $Ext^1_{U}(Q, M)$.
Therefore $O_a(U, M)$ is finitely generated over $A$ as well
and by the same corollary it is independent of
$U$ as long as $U \cap Z^\circ_A = \emptyset$. 
For such $U$, denote the stabilized module
$O_a(U, M)$ simply by $O_a(M)$.

\bigskip
\noindent
Assume now that $A$ is itself an extension
$$
0 \to N \to A \to A_0 \to 0
$$
where $N$ is nilpotent and acts on $M$ by zero (so that $M$ is an $A_0$-module).
Both $Ext^1_{U}(Q, M)$ and $Ext^1_{U}(L_A, M)$ in this case are 
$A_0$-modules,
hence $O_a(M)$ is also a finitely generated $A_0$-module. In fact, if 
$a_0 = (E_0, U_0)$ is the restriction of $a = (E, U)$ and 
$Q^\bullet$ and $L^\bullet$ are the above locally free resolutions, then
$$
Ext^1_{U}(Q, M)= Ext^1_{U_0}(Q^\bullet \otimes_A A_0, M); \qquad 
Ext^1_{U}(L(A), M) = Ext^1_{U_0}(L^\bullet \otimes_A A_0, M)
$$
and by virtue of the spectral sequence just mentioned we can assume
that $Q^\bullet$, $L^\bullet$ are concentrated in degrees $[-2, 0]$.
Denoting $P^\bullet = \mathcal{H}om_{U_0}(Q^\bullet \otimes_A A_0, \mathcal{O}_{U_0}) = (P^0 \to P^1 \to P^2)$ and similarly for 
$R^\bullet = \mathcal{H}om_{U_0}(L^\bullet \otimes_A A_0, \mathcal{O}_{U_0})
= (R^0 \to R^1 \to R^2)$ we see that 
$$
O_a(M)= Coker \Big[H^1(U_0, P^\bullet_M) \to H^1(U_0, R^\bullet_M)\Big]
$$
By the standard argument, cf. [I], if the obstruction $\omega(E)$
vanishes, all deformations of $E$ over $U'$ form a pseudo-torsor 
over $D_{a_0}(M)$.

\subsection{Etale localization, completions, constructibility.}

Let $p: A \to B$ be etale and consider $p_0: A_0 \to B_0$ defined by 
$B_0 = B\otimes_A A_0$, $p_0 = p \otimes_A A_0$. Consider 
$a = (E, U) \in F(A)$ and let $b$ be its pullback in $F(B)$, and 
similarly for $a_0 = (E_0, U_0) \in F(A_0)$, $b_0 \in F(B_0)$. 
\begin{prop} There exist natural isomorphisms:
$$
O_b(M\otimes_{A_0} B_0) \simeq O_a(M)\otimes_{A_0} B
$$
and
$$
Aut_{b_0}(B_0 + M\otimes_{A_0} B_0) \simeq 
Aut_{a_0}(A_0 + M) \otimes_{A_0} B_0; \quad D_{a_0} (M\otimes_{A_0} B_0)
\simeq D_{a_0} (M)\otimes_{A_0} B_0
$$
\end{prop}
\textit{Proof.} For $Aut$ and $D$ this
follows from their identification with 
$H^i(U_0, End(E_0)_M)$
for $i=0, 1$ and some $U_0 \subset X_{A_0}$, and etale localization 
for cohomology.
For $O$, one uses the definition
$$
Coker\Big[Ext^1_{U}(Q, M) \to Ext^1_{U}(L(A), M)\Big]
$$
with some $U\subset X_{A}$, and then applies etale localization for
$Ext^1$ plus the idenity $L(B) = L(A) \otimes_A B$ which holds for any 
etale extension $A\to B$, cf. Chapter II of \cite{I}. $\square$

\bigskip

\begin{prop}
Let $\mathfrak{m} \subset A_0$ be a maximal ideal. Then
$$
D_{a_0}(M) \otimes_{A_0} \widehat{A}_0 \simeq 
\underleftarrow{\lim}\; D_{a_0}(M/ \mathfrak{m}^nM)
$$
and similarly for $Aut_{a_0} (A_0 + M)$. 
\end{prop}
\textit{Proof.} Both follow
immediately from Proposition 0.13.3.1 in [EGA $\textrm{III}_1$] applied
to the completion of the open subscheme $U_0$ for which 
$H^i(U_0, End(E_0)_M)$ compute for $i = 0, 1$ the modules 
$Aut_{a_0}$ and $D_{a_0}$, respectively. $\square$

\bigskip

\begin{prop}
Assume that the ring $A_0$ is \textrm{reduced}. Then there 
exists an open dense subset of points of finite type $p \in Spec(A_0)$,
so that
$$
D_{a_0}(M)\otimes_{A_0} k(p) \simeq
D_{a_0}(M \otimes_{A_0} k(p)),
$$
and similarly for $Aut_{a_0}(A_0 + M)$ and $O_a(M)$.
\end{prop}
\textit{Proof.} 
\textit{Step 1.}
Without loss of generality we can also assume that
$Spec(A_0)$ is irreducible, i.e. $A_0$ is a domain. 
It suffices to show
that localizing $Spec(A_0)$ at 
(powers of) a nonzero element one can achieve
$$
H^i(U_0, P^\bullet_M) \simeq H^i(U_0, P^\bullet) \otimes_{A_0} M; \quad
i = 0, 1
$$
for arbitrary finitely generated $A_0$-module $M$ and a fixed complex
$P^\bullet = (P^0 \to P^1 \to P^2)$ of vector bundles on $U_0$.

First consider the case when $P^\bullet = P^0 =:P$ is a 
vector bundle in degree zero. To unload notation write $X_0$ for $X_{A_0}$, 
denote the open 
embedding $U_0 \to X_0$ by $j$ and recall that $j_* (P_M)$, 
$R^1 j_*(P_M)$ are coherent by Corollary \ref{finitely}. 
For a quasicoherent sheaf $\mathcal{F}$ on $X_0$
recall also the exact sequence:
$$
\ldots \to H^i_Z(X_0, \mathcal{F}) \to H^i(X_0, \mathcal{F}) \to H^i(U_0, j^* \mathcal{F}) \to H^{i+1}_Z(X_0, \mathcal{F}) \to \ldots
$$
and the spectral sequence $E_2^{p, q} = H^q(X_0, \mathcal{H}^p_Z(\mathcal{F}))
\Rightarrow H^{p+q}_Z(X_0, \mathcal{F})$. Since $\mathcal{F} = j_*(P_M)$
satisfies $\mathcal{F} \simeq j_* j^* \mathcal{F}$ 
we have $\mathcal{H}^0_Z(j_*(P_M)) = \mathcal{H}^1_Z(j_* (P_M)) = 0$.  
This gives an isomorphism
$
H^0(U_0, P_M) \simeq H^0(X_0, j_*(P_M))
$
and a long exact sequence 
\begin{equation}
\label{long}
0 \to H^1(X_0, j_*(P_M)) \to H^1(U_0, P_M) 
\to H^0(X_0, R^1j_*(P_M)) \to H^2(X_0, j_*(P_M))
\end{equation}

\bigskip
\noindent
\textit{Step 2.}
Since $X_0$ is proper over $Spec(A_0)$, for finite complex of vector bundles 
$\mathcal{F}$ on $X_0$ and any finitely generated $A_0$-module $M$ one has
\begin{equation}
\label{onX}
H^i(X_0, \mathcal{F}_M) 
\simeq H^i(X_0, \mathcal{F}) \otimes_{A_0} M, \quad i \geq 0
\end{equation}
after a localization of $A_0$ at a nonzero element $f \in A_0$. 
In fact, only finitely many cohomology modules are non-zero 
and by properness they are finitely generated over $A_0$. Applying the following
Generic Freeness Lemma, cf. Theorem 14.4 in \cite{E}, we
can ensure that after a localization
of $A_0$ the  cohomology modules will be free of finite rank:
\begin{lemma}
Let $A_0$ be a noetherian domain and $B$ a finitely generated $A_0$-algebra.
If $N$ is a finitely generated $B$-module, there exists a nonzero element
$t \in A_0$, such that the localization $N[t^{-1}]$ is free over $A_0[t^{-1}]$.
$\square$
\end{lemma}
Assuming that all $H^i(X_0, \mathcal{F})$ are free of finite rank over $A_0$, 
let   $C^\bullet(\mathcal{F})$ be the Cech complex of $\mathcal{F}$.
It is a complex of flat 
$A_0$-modules, such that $C^\bullet(\mathcal{F}) \otimes_{A_0} M$ computes the 
cohomology 
of $\mathcal{F}_M$ for any $M$. This gives a second quadrant spectral 
sequence
$$
E_2^{p, -q} = Tor^{A_0}_q(H^p(X_0, \mathcal{F}), M) \Rightarrow 
H^{p-q}(X_0, \mathcal{F}_M)
$$
which by freeness of $H^i(X_0, \mathcal{F})$ reduces to 
$H^i(X_0, \mathcal{F}_M) \simeq
H^i(X_0, \mathcal{F}) \otimes_{A_0} M$, as required.

\bigskip
\noindent
\textit{Step 3.} Denote the exact sequence \eqref{long} by $K^\bullet(M)$.
Localizing $A_0$  we can assume that all modules in $K(A_0)$ and 
$$
Coker \big[H^0(X, R^1j_*P) \to H^2(X, j_*P)\big]
$$
are free over $A_0$. Then  $K^\bullet(A_0) \otimes_{A_0} M$ is again exact.
Comparing $K^\bullet(A_0)\otimes_{A_0} M$ with $K^\bullet(M)$ and 
using the isomorpism of Step 2, we reduce the isomorphisms 
\begin{equation}
\label{change}
H^i(U_0, P_M) \simeq H^i(U_0, P)\otimes_{A_0} M \quad i = 0, 1
\end{equation}
to their local versions:
\begin{lemma} In the notation introduced above, 
\begin{equation}
\label{local}
j_*(P_M)\simeq j_*(P)_M, \qquad R^1j_*(P_M) \simeq R^1j_*(P)_M
\end{equation}
\end{lemma}
\textit{Proof of the lemma.} Since the statement in local we can 
assume that $X_0 = Spec(B)$ is affine.
Let $N$ be the $B$-module corresponding to 
the coherent sheaf $j_*(P)$ and $I \subset B$ an ideal such that 
$Supp(B/I) = Z$. Then the local cohomology modules $H^0_I(N)$, $H^1_I(N)$ 
vanish by Step 1. Hence by Proposition 18.4 in \cite{E}, 
$depth_I(N) \geq 2$ and there
is an $N$-regular sequence $(f, g) \in I$. We can also
assume that $(f, g)$ is regular on $N_M = N\otimes_{A_0} M$  for any $M$.
 Indeed, by $N$-regularity we have an exact sequence
$$
0 \to N \to N \oplus N \to N \to N/(f, g)N \to 0
$$
of finitely generated $B$-modules. Localizing $A_0$ we may assume that all
 modules in this sequence are free over $A_0$ therefore the complex 
$$
0 \to N_M \to N_M \oplus N_M \to N_M \to 
\big(N/(f, g)N\big)\otimes_{A_0} M \to 0
$$
is also exact. Since its first three terms give the Koszul complex of $N_M$,
$(f, g)$ is $N_M$-regular. 
In particular, $H^0_I(N_M) = H^1_I(N_M) = 0$. Since $N_M$ corresponds to the
sheaf $j_*(P)_M$, the vanishing of local cohomology gives 
$$
j_*(P)_M \simeq j_* j^*(j_*(P)_M) 
\simeq j_*(P_M);
$$
i.e. the first isomorphism of \eqref{local}.
The second isomorphism  is equivalent to 
$$
H^2_I(N_M) \simeq H^2_I(N) \otimes_{A_0} M.
$$
We claim that one can replace $I$ by an ideal $I' = (f, g, h) \subset I$, 
where $h \in I$, such that $H^2_I(N_M) \simeq H^2_{I'}(N_M)$ for all $M$.
In fact, by prime avoidance, cf. Lemma 3.3 in \cite{E}, one can find $h \in I$ 
such that $h \notin Q$ whenever $Q$ is an associated prime of $B/(f, g)B$ 
not containing $I$. Since $P$ is locally free away from $Z$, the 
sequence $(f, g)$ is $\mathcal{O}_U$-regular.
This implies that $Z' = V(f, g, h)$ has codimension 3 in $X$ at any 
point of of $Z' \setminus Z$. Recall that by our assumption 
$Z$ contains all points where the fiberwise $S_3$ conditon is violated.
Therefore by Lemma \ref{serre} 
the sheaf $j_*(P)_M$ corresponding to $N_M$ satisies the relative 
$S_3$ condition on $U$ and 
$\mathcal{H}^2_{Z'}(j_*(P)_M)|_{U} = 0$. 
Then the sheaf $\mathcal{H}^2_{Z'}(j_*(P)_M)$ is supported at $Z$ and
$$
H^2_{I'}(N_M) \simeq H^0_I(H^2_{I'}(N_M)) \simeq H^2_I(N_M)
$$
where the last isomorphism uses the spectral sequence $H^p_I(H^q_{I'}(N_M))
\Rightarrow H^{p+q}_I(N_M)$ coming from $R\Gamma_{I_1 + I_2} \simeq
R\Gamma_{I_1} \circ R\Gamma_{I_2}$; and vanishing of $H^0_{I'}(N_M)$ and
$H^1_{I'}(N_M)$ due to  $N_M$-regular sequence $(f, g)\in I'$. 

The local cohomology $H^i_{I'}(N_M)$ for $i \geq 2$
are computed by the complex
$C^\bullet(P) \otimes_{A_0} M$ where $C^\bullet(P)$ is the Cech complex
of the vector bundle $P$ on $X \setminus Z'$ with respect to the affine covering 
$X_f \cup X_g \cup X_h$. Observe that  each term of $C^\bullet(P)$ is flat 
over $A_0$, and the tensor product $C^\bullet(P)\otimes_{A_0} M$ can be 
identified with the Cech complex of $P_M$. Then the second quadrant
spectral sequence 
$$
E_2^{-p, q} = Tor_p^{A_0}(H^q(C^\bullet(P)), M) \Rightarrow H^i(C^\bullet(P)
\otimes_{A_0} M)
$$
reduces to $H^q_{I'}(N_M)\simeq H^q_{I'}(N) \otimes_{A_0} M$  for
any $M$ and $q\geq 2$ 
if all $H^q_{I'}(N)$ are flat over $A_0$ (recall
that  for $q = 0,1$ and
these modules vanish by the above discussion). For $q=2$ we have
$$
H^2_{I'}(N) = H^2_I(N) = \Gamma(X, R^1j_* P)
$$
which is finitely generated over $B$ since $R^1j_*P$ is coherent. By 
Generic Freeness we can assume that $H^2_{I'}(N)$ is free over $A_0$
after a localization. Since the Cech complex only has three term, 
the only remaining local cohomology group is
$$
H^3_{I'}(N) = \underrightarrow{\lim} \; N/(f^n, g^n, h^n) N
$$
and this is not finitely generated over $B$. 
Since a 
filtered direct limit of projectives is flat, 
it suffices to show that all $N_n = N/(f^n, g^n, h^n)N$,  $n \geq 1$, 
become projective after a single localization of $A_0$.  

By a standard combinatorial argument $N_n$ admits a filtration with 
quotients of the type 
$$
N_{p, q, r}=
f^p g^q h^r N/(f^{p+1} g^q h^r, f^p g^{q+1} h^r, f^p g^q h^{r+1}) N.
$$
and it suffices to ensure projectivity of these modules since $N_n$
are their iterated extensions.
Both $B$ and $N$ have triple filtrations by powers of $f, g, h$ 
respectively and taking associated graded objects $gr(\cdot) =
gr_h gr_g gr_f(\cdot)$
we obtain a $\mathbb{Z}^3$-graded ring $gr(B)$ and a finitely generated
module $gr(N) = \bigoplus_{p, q, r \geq 0} N_{p, q, r}$. Applying 
Generic Freeness to $gr(N)$ we can localize $A_0$ at a single element
and assume that $gr(N)$ is free over $A_0$. Then each direct factor
$N_{p, q, r}$ of $gr(N)$ must be projective over $A_0$, as required.
This proves the lemma. $\square$

\medskip
\noindent
Therefore, for a vector bundle $P$ on $U$ we have established isomorphisms
$$
H^i(U_0, P_M) \simeq H^i(U_0, P) \otimes_{A_0} M; \qquad
i = 0, 1.
$$

\bigskip
\noindent
\textit{Step 4.}
Now we return to the general case of 
a complex $P^\bullet = (P^0 \to P^1 \to P^2)$
of vector bundles on $U_0$. Consider the spectral sequence 
$
E_1^{r, q}(M)= H^q(U_0, P^r_M) \Rightarrow H^{r+q}(U_0, P^\bullet_M).
$
Then $E_2^{r, q}$ is the $r$-th cohomology of 
$$
E^{\bullet, q}_1 (M) = \big[H^q(U_0, P^0_M) \to H^q(U_0, P^1_M) \to 
H^q(U_0, P^2_M)\big]
$$
By previous step for $q = 0, 1$ we can localize $A_0$ to achieve
$H^q(U_0, P^r_M) \simeq H^q(U_0, P^r) \otimes_{A_0} M$. In addition, we 
can ensure that the cohomology of the complexes $E^{\bullet, q}_1(A_0)$ for
$q = 0, 1$, are free finitely generated $A_0$-modules. The second assumption
guarantees that for $E^{\bullet, q}_1(A_0)$ cohomology commutes with 
$\otimes_{A_0} M$; which in view of the first assumption gives $E^{r, q}_2(M) 
= E^{r, q}_2(A_0) \otimes_{A_0} M$ for $q = 0, 1$.
Now the proposition follows by the isomorphism 
$H^0(U_0, P^\bullet_M) = E^{0, 0}_2(M)$ and the exact sequence 
$$
0 \to E^{1, 0}_2(M) \to H^1(U_0, P^\bullet_M)\to E^{0, 1}_2(M) \to E^{2, 0}_2(M)
\hspace{2cm} \square
$$

\section{Effectiveness}

\begin{prop}
Let $\widehat{A}$ be a complete local algebra with residue field of finite
type over $k$ and maximal ideal $\mathfrak{m}$, then the canonical functor
$$
F(\widehat{A}) \to \underleftarrow{\lim}\;F(\widehat{A}/\mathfrak{m}^n) 
$$
is an equivalence of categories.
\end{prop}
\textit{Proof.} Let $\{(E_n, U_n)\in F(\widehat{A}/\mathfrak{m}^n)\}$ be 
a sequence representing an object on the right hand side. 
Shrinking $U_n$ as in Section 2 we can assume that each $E_n$
satisfies the relative $S_3$ condition on $U_n$. But then by the
tangent-obstruction theory  and stabilization of cohomology
the set of isomorphism classes of extensions of $E_n$
to $F(A/\mathfrak{m}^{n+1})$, i.e. the cohomology 
$H^1(?, End(E_n) \otimes \mathfrak{m}^n/\mathfrak{m}^{n+1})$,
will be the same over $U_n$ as over any open subset 
$W \subset U_n \cap U_{n+1}$ with closed complement in 
$\Phi(\widehat{A}/\mathfrak{m}^n)$.
Therefore we can assume that $E_{n+1}$ is defined also on $U_n$ 
and by induction all $E_n$ are defined on the
same open subset $U'$. Then by the main result of [B] there
exists a bundle $E$ on an open subset $U \subset Spec(\widehat{A})$
such that $(E, U)$ restricts to $(E_n, U')$ in each 
$F(A/\mathfrak{m}^{n})$. On morphisms the assertion also follows from 
the main result of \textit{loc.cit}. $\square$

\section{Properties of the diagonal}

\begin{lemma}
\label{q}
Let $G$ be a coherent sheaf on $X_k$
and $H$ a coherent sheaf on $X_A$. Then there exists a finitely 
generated $A$-module $Q$, unique up to canonical isomorphism, 
and a natural isomorphism of covariant functors (with argument M)
$$
Hom_{X_A} (H, (G \otimes_k A) \otimes_A M) \simeq Hom_A(Q, M)
$$
from the category of $A$-modules to itself. 
\end{lemma}
\textit{Proof.} First assume that $X_k$ is projective. Since
$G \otimes_k A$ is flat over $A$ and $H$ is a cokernel of a 
morphism of locally free sheaves, the assertion is an immediate
consequence of Corollary 7.7.8 in [EGA $\textrm{III}_2$].

For proper $X_k$ we use Chow Lemma and the pattern
of Section 5 in [EGA $\textrm{III}_1$]. Fixing $H$, we
call $G$ \textit{representable} if a module $Q$ as in the statement
exists. Consider an exact sequence if coherent sheaves on $X_k$:
$$
0 \to G_1 \to G_2 \to G_3 \to 0.
$$
We claim that if $G_2$, $G_3$ are representable then $G_1$ is, 
and if $G_1, G_3$ are representable then $G_2$ is. The first 
assertion is quite easy as the morphism $G_2 \to G_3$ 
corresponds to a morphism $Q_3 \to Q_2$ of representing 
$A$-modules and hence $Q_1 = Coker(Q_3 \to Q_2)$ represents $G_1$. For the
second assertion we first show that the functor
$M\mapsto R(M)= Hom_{X_A} (H, (G \otimes_k A) \otimes_A M)$ commutes
with projective limits. In fact, choose an affine covering 
$X_k = \cup U_i$ and compute $R(M)$ as the kernel of the 
first arrow in the corresponding Cech complex
$$
\bigoplus_i Hom_{(U_i)_A}(H, (G \otimes_k A) \otimes_A M)
\to 
\bigoplus_{i \neq j} 
Hom_{(U_i \cap U_j)_A}(H, (G \otimes_k A) \otimes_A M)
$$
where we suppress the notation for restriction of sheaves to $U_i$
and $U_i \cap U_j$, respectively. Since projective limits
are left exact and commute with $Hom(H, \cdot)$ by universal 
property of projective limits, it suffices to show that
$(G\otimes_k A) \otimes_A(\cdot)$ commutes with projective
limits, which is obvious since the first tensor factor is 
free over $A$. By a theorem of Watts, \cite{W}, since $R(M)$ is
left exact and commutes with projective limits, it is 
representable by an $A$-module $Q$: $R(M) = Hom_A (Q, M)$ 
although in general $Q$ may not finitely generated. But the 
exact sequence of $G_i$ induces an exact sequence
of representing modules
$$
Q_1 \to Q_2 \to Q_3 \to 0
$$
and if $Q_1$, $Q_3$ are finitely generated, the same 
holds for $Q_2$. 

Now we prove the assertion for all proper separated $X_k$ 
by induction on the dimension of $Supp(G)$ (as before, 
$A$ and $H$ are fixed). Passing to an appropriate closed
subscheme in $X_k$ we can assume $Supp(G)= X_k$. By 
Chow Lemma there exists a projective 
scheme $\widetilde{X}_k$ over $k$ and a projective morpism
$h: \widetilde{X}_k \to X_k$ which is an isomorphism 
over a dense open subset of $X_k$. Let $h_A: \widetilde{X}_A 
\to X_A$ be the morphism obtained by base 
change $Spec(A)\to Spec(k)$. By adjunction
$$
Hom_{\widetilde{X}_A} (h_A^* H, (h^*(G) \otimes_k A) \otimes M)
\simeq Hom_{X_A}(H, (h_* h^*(G) \otimes_k A) \otimes M)
$$
so the coherent sheaf $h_* h^*(G)$ is representable. Since
the kernel and the cokernel of $\phi: G \to h_* h^* (G)$ are
zero on a dense open subset of $X_k$, by induction they 
are representable. Hence by above argument $Im(\phi) = 
Ker(h_* h^*(G) \to Coker(\phi))$ is representable and 
from the exact sequence $0 \to Ker(\phi) \to G \to Im(\phi) \to 0$
the sheaf $G$ is also representable, as required. $\square$
\begin{cor}
\label{sec}
Let $E$ be a vector bundle on $U= X_A \setminus Z$ with 
$Z\in \Phi(A)$ and $Y \to U$ a closed subscheme in the total 
space of $E$ over $U$. 

(i) The functor $Sec(E/U)$ on 
$(Aff/A)$ which sends an $A$-algebra $B$ to the set of all
sections $U_B \to E_B$ is 
represented by an affine scheme $Spec(Sym^\bullet_A(Q))$
for a finitely generated $A$-module $Q$.

(ii) The similar functor $Sec(Y/U)$ is represented by a
closed subscheme of $Spec(Sym^\bullet_A(Q))$. 
\end{cor}
\textit{Proof.} First we deal with $Sec(E/U)$. By definition
$$
Sec(E/U)(B)= 
Hom_{\mathcal{O}_U-alg}
(Sym^\bullet(E^\vee), \mathcal{O}_U  \otimes_A B) 
= Hom_{\mathcal{O}_U}
(E^\vee, \mathcal{O}_U  \otimes_A B) 
$$
If $j: U\to X_A$ is the open embedding, then combining the adjunction
of $j_*$, $j^*$  with  $j^*j_*(E^\vee) \simeq E^\vee$, $j_*(\mathcal{O}_U
\otimes_A B) \simeq \mathcal{O}_{X_A} \otimes_A B$ we get
$$
Hom_{\mathcal{O}_U}
(E^\vee, \mathcal{O}_U  \otimes_A B) 
= Hom_{\mathcal{O}_{X_A}}
(j_*(E^\vee), \mathcal{O}_{X_A}\otimes_A B).
$$
By the previous lemma, for some finitely generated $A$-module $Q$
we can identify the last module with 
$$
Hom_{A}(Q, B) = Hom_{A-alg}(Sym^\bullet_A(Q), B);
$$
hence the functor $Sec(E/U)$ is represented by $Spec(Sym^\bullet_A(Q))$.

For a closed subscheme $Y$ we can find a $\mathcal{O}_U$-coherent
subsheaf $N$ of $Sym^\bullet_{\mathcal{O}_U} (E^\vee)$ which 
generates the ideal subsheaf of $Y$ as 
$Sym^\bullet_{\mathcal{O}_U} (E^\vee)$-module, e.g. by following
the pattern of Proposition 9.6.5 in [EGA I]. 
A section $s: U_B \to E_B$ induces a homomorphism of 
$\mathcal{O}_U$-algebras $Sym^\bullet_{\mathcal{O}_U}(E^\vee)
\to \mathcal{O}_U\otimes_A B$ and $s$ factors through $Y_B$
precisely when the restriction 
$\phi: N \to \mathcal{O}_U \otimes_A B$ vanishes.
Fixing a coherent sheaf $N'$ on $X_A$ with $j^*(N')\simeq N$
and using the adjunction one more time, 
we get an isomorphism
$$
Hom_{X_A}(N', \mathcal{O}_{X_A} \otimes B) \simeq
Hom_{U}(N, \mathcal{O}_{U} \otimes_A B)
$$
By Lemma \ref{q} there exists a finitely generated $A$-module $R$
such that the above $Hom$ groups can be identified with $Hom_A(R, B)$. 
Denote the corresponding homomorphism $R \to B$ by the same letter $\phi$. 
Then for any $B$-algebra $B\to B'$ the induced section 
$s_{B'}$ corresponds to the composition of $\phi:R \to B$
with $B \to B'$. It follows that $s_{B'}$ factors through $Y_{B'}$ 
presicely when $Ker(B\to B')$ contains the ideal generated by 
$\phi(R) \subset B$. Therefore $Sec(Y/U)$ is a closed subfunctor 
of $Sec(E/U)$ and the assertion follows.  $\square$.

\begin{prop}
The diagonal of $F$ is representable, quasi-compact and separated. 
\end{prop}
\textit{Proof.} Although representability of the diagonal follows 
formally from the previous results, it is useful to establish it directly:
if $S$ is an algebraic space then a morphism $S \to F \times_k F$ 
corresponds to a pair of rank $r$ bundles $E_1, E_2$ which we may assume to 
be defined on a common open subset $X_S \setminus Z$. Then the 
fiber product with the diagonal is the functor of isomorphisms  
$Isom(E_1, E_2)$. Although such isomorphisms correspond to sections with 
values in an \textit{open} subset of the vector bundle $Hom(E_1, E_2)$, 
the isomorphism condition is equivalent to the nonvanishing of
the determinant, i.e. the induced section of 
$L = Hom(\Lambda^rE_1, \Lambda^r E_2)$. Therefore the subfunctor of 
isomorphisms can be identified with the closed subscheme in the total space
of $Hom(E_1, E_2) \oplus L^*$, formed by all sections $(\phi, s)$ 
such that $\det(\phi) s = 1$. Now Corollary \ref{sec} gives 
representability in the case of affine $S$, and uniqueness of
the representing module $Q$ from Lemma \ref{q}, in the general case. 
Quasi-compactness follows immediately from the fact that $Isom(E_1, E_2)$
is affine of finite type over $S$. 

By valuative 
criterion, cf. Theorem 7.3 in \cite{La}, separatedness reduces to the following
fact: if $R$ is a discrete
valuation ring and $(E_1, U_1)$, $(E_2, U_2)$ two objects in $F(R)$
then $H^0(U_1 \cap U_2, \mathcal{H}om (E_1, E_2))$ is torsion free.
Denote $U = U_1 \cap U_2$, $E= Hom(E_1, E_2)$ and let 
 $t \in R$ be a local parameter.
Then we need to show that $ts = 0$ for $s \in H^0(U, E)$ implies $s=0$.
This question is local on $U$ hence we can assume that $E$ is a trivial 
bundle and $U= Spec(B)$, for a flat $R$-algebra $B$. Then a short exact
sequence $0 \to R \stackrel{t}{\longrightarrow} R \to K \to 0$ gives
$0 \to B  \stackrel{t}{\longrightarrow} B \to B \otimes_R K \to 0$,
which proves the assertion. $\square$

\bigskip
\noindent
\textbf{Remark.} Observe that as in
Corollary \ref{sec}, for any pair of bundles $E, F$ on $U= X_A \setminus Z$
the functor on $A$-modules, which sends $M$ to 
$Hom_U(E, F \otimes_A M)$ is represented by a finitely generated $A$-module
$Q$ which is free of finite rank over a dense open subset of $Spec(A)$. 
For a general coherent sheaf $M$ on an algebraic space $S$ we can glue
the representing modules using uniqueness of $Q$ and obtain a vector 
bundle $Q$ on a dense open subset $S_0 \subset S$, which represents $Hom(E, F \otimes_{\mathcal{O}_S}M)$ over $S_0$. This is a
direct analogue of  Proposition 2.2.3(i) in \cite{L}.


\section{Representability of principal bundles.}

\textit{Proof of Theorem 1.} To prove that $F_{GL(r)}$ is an algebraic stack,
locally of finite type and separated over $k$, we just need to collect 
the results of the previous sections and compare with the conditions
of Artin's representability criterion, cf. Theorem 5.3 in \cite{Ar}. The ``limit 
preserving" condition is proved in Section 3, Schlessinger's condition S1
in Section 4, while S2 is establihsed in Section 5.1. Effectiveness (condition
(2) of Artin's criterion) is given by Section 3, while part (3) 
of \textit{loc.cit.} is proved in Section 5.2. Finally local quasi-separation
(part (4) of Artin's criterion) is established in a stronger form in Section 7.

For a general reductive group $G$ over $k$ we use a result of Haboush, 
cf. \cite{Ha},
and choose an exact finite dimensional representation
$\rho: G \to GL(r)$ with $Y = GL(r)/G$ affine. Moreover, $Y$ is
isomorphic to a closed $GL(r)$-orbit of a vector in a finite dimensional 
rational $GL(r)$-module $W$.

Then each principal $G$-bundle $P$
induces a principal $GL(r)$-bundle $E = P_\rho$. Conversely, for any principal 
$GL(r)$-bundle $E$ over a scheme $U$ 
its reduction of the structure group to $G$ may be
viewed as a regular section $U \to Y_E = E \times_{GL(r)} Y$. 
Moreover, $Y_E$ is a closed subscheme of a total space of a vector bundle
$W_E$ on $U$, induced from $E$ via the homomorphism $GL(r) \to GL(W)$. 

This construction and Corollary \ref{sec} shows that the morphism $F_G \to F_{GL(r)}$ is
representable in the sense of Definition 3.9 in \cite{La} and since
$F_{GL(r)}$ is an algebraic stack, by Proposition 4.5 (ii) of \textit{loc.cit.}
$F(G)$ is also an algebraic stack.

Alternatively, we can re-prove the results of Sections 3, 4 and 7 for $F_G$ by 
using the proved facts for vector bundles and reducing the 
structure group from $GL(r)$ to $G$. The Effectiveness property of Section 6
is proved in \cite{B} by a similar strategy. Finally, by Chapter VI of
\cite{I} the arguments of Section 5 carry over to $G$ after a minor 
modification. Let $\Omega(G)$ be the space of $G$-invariant (from the right) 
differential forms on $G$, with the natural left $G$-action 
(in characetistic zero this is just the adjoint representation). For 
any $G$-bundle $P$ denote by $ad(P)$ the vector bundle induced from 
$P$ via the action homomorphism $G \to GL(\Omega(G))$. Then all arguments 
of Section 5 carry over if $Ext^i(End(E), \cdot)$ are 
replaced by $Ext^i(ad(P), \cdot)$ and $H^i(\cdot, End(E)_M)$ by 
$H^i(\cdot, ad(P)^\vee_M)$, where $ad(P)^\vee$ stands for the dual 
bundle. This finishes the proof of Theorem 1. $\square$

\bigskip
\noindent
\textbf{Remark.} If we replace the ``codimension 3" condition in the
definition of $\Phi(T)$ (see Section 1), by ``codimension 2", the stack
$F_G$ will not longer be algebraic. The most obvious reason is that 
the tangent space $H^1_{A_0, \Phi}(End(E))$ will no longer be finitely
generated over $A_0$. On the other hand, the results of Sections 3, 4 and 7
remain valid while the Effectiveness of Section 6, which definitely fails
as such, may be repaired by introducing an additional condition as 
 in \cite{B}. It is conjectured by V. Drinfeld that in the codimension 2 case
$F_G$ is an inductive limit of algebraic stacks, locally of finite type over
$k$. We plan to return to this topic, as well as the related Uhlenbeck functor, 
in future work.

\bigskip
\noindent
\textit{Address:} Department of Mathematics, UC Irvine, Irvine CA, 92697.

\noindent
\textit{Email:} vbaranov@math.uci.edu

\end{document}